\documentclass[11pt]{article}
\usepackage{fullpage}
\usepackage{amsthm}
\usepackage{comment}
\usepackage{amsmath}
\usepackage{amssymb}
\usepackage{algorithm}
\usepackage{algorithmic}
\usepackage{bm}
\usepackage{comment}
\usepackage{hyperref}
\usepackage{systeme}

\newtheorem{theorem}{Theorem}[section]
\newtheorem{lemma}{Lemma}[section]

\newtheorem{proposition}{Proposition}[section]

\newtheorem{conj}{Conjection}
\newtheorem{defi}{Definition}[section]

\title{{\small Algebraic Geometry, Number Theory, Moufang loops}\\
An example of a non-associative Moufang loop of point classes on a cubic surface}
\author{D. Kanevsky\\ 
}

\begin{document}

\newcommand{\Addresses}{{
  \bigskip
  \footnotesize

  D. Kanevsky, \textsc{Google Research, 1600 Amphitheatre Parkway, Mountain View, CA 94043}\par\nopagebreak
  \textit{E-mail address}, \texttt{dkanevsky@google.com}

}}
\maketitle
\Addresses
\begin{center}
    \large
    \textbf{Abstract}
\end{center}
Let $V$ be a cubic surface defined by the equation $T_0^3+T_1^3+T_2^3+\theta T_3^3=0$ over a quadratic extension of 3-adic numbers $k=\mathbb{Q}_3(\theta)$, where $\theta^3=1$. We show that a relation on a set of geometric k-points on $V$ modulo  $(1-\theta)^3$ (in a ring of integers of $k$) defines an admissible relation and a commutative Moufang loop  associated with classes of this admissible equivalence is non-associative. This answers a problem that was formulated by Yu. I. Manin more than 50 years ago about existence of a cubic surface with a non-associative Moufang loop of point classes. 
\section{Introduction}
\numberwithin{equation}{section}
\label{intro}
Throughout the paper $V$ will denote an irreducible cubic hypersurface defined over a field $k$ and embedded in some projective space  over $k$. 
It is well-known that new points on $V$ may be obtained from old ones by secant and tangent constructions. One of difficulties in studying this secant and tangent process when $\dim V > 1$ is that it does not lead to an associative binary operation as in the case of elliptic curves. An interesting remedy to this was suggested in (\cite{siksek}) via introduction of a certain abelian group associated with a cubic surface which retains much information about secant and tangent process. Here we follow a different approach that was suggested in (\cite{manin3}). See also a short survey of  content of first two chapters of this Manin's book in Section \ref{summary} in this paper.\\ 
Let $V(k)$ be a set of points on $V$ (where by ``points" we mean geometric points with values in $k$). 
 Three points $\bm{P_1}, \bm{P_2}, \bm{P_3} \in V(k)$ not necessary different, are said to be collinear if they lie on a straight line $\bm{L}$ defined over $k$ and either $\bm{L} \subset V$ or $\bm{P_1}+\bm{P_2}+\bm{P_3}$ is the intersection cycle of $V$ with $\bm{L}$. For any $\bm{P_1}, \bm{P_2} \in V(k)$ there exists $\bm{P_3} \in V(k)$ such that $\bm{P_1},\bm{P_2}, \bm{P_3}$ are collinear.  
 If $V$ is smooth and has dimension 1 then the collinearity relation on its set of points introduces  a structure of the Abelian group by the law $\bm{P}_1+\bm{P}_2 = \bm{U} \circ (\bm{P}_1 \circ \bm{P}_2)$ where $\{\bm{P}_1, \bm{P}_2, \bm{P}_1 \circ \bm{P}_2\}$ and $\{\bm{P}_1+\bm{P}_2, \bm{P}_1 \circ \bm{P}_2, \bm{U}\}$ are collinear triples of points on $V$).
If $V$ has dimension $ \geq 2$ than the condition ``$\bm{P_1}, \bm{P_2}, \bm{P_3} \in V(k)$ are collinear'' does not define a binary composition since there are many tangent lines at a point $\bm{P_1}$ to $V$ or there could be lines through $\bm{P_1}$ belonging to $V$. 
This leads to the following definition.
\begin{defi}
\label{admDef} 
An equivalence relation $\mathcal{A}$ on $V(k)$ is said to be admissible if the following condition holds:\\
If $\bm{P_1}, \bm{P_2}, \bm{P_3} \in V(k)$ and $\bm{P_1'}, \bm{P_2'}, \bm{P_3'} \in V(k)$ are collinear, $\bm{P_1} \sim \bm{P_1'} \mod \mathcal{A}$,
$\bm{P_2} \sim \bm{P_2'} \mod \mathcal{A}$ then $\bm{P_3} \sim \bm{P_3'} \mod \mathcal{A}$.
\end{defi}
 If $\bm{P_1} \neq \bm{P_2}$ and there is no straight line on $V$ through $\bm{P_1}, \bm{P_2}$ then $\bm{P_3}$ is uniquely defined and we sometimes write  $\bm{P_1} \circ \bm{P_2} = \bm{P_3}$. 
Let $\mathcal{A}$ be an admissible equivalence on $V(k)$ and let $E = V(k)/\mathcal{A}$ denote the set of equivalence classes. For any $\bm{S_1}, \bm{S_2} \in E$ there is a uniquely determined class $\bm{S_3} \in E$ such that that there exists collinear points $\bm{P_1}, \bm{P_2}, \bm{P_3}$ with $\bm{P_i} \in \bm{S_i}$. Defining a new binary operation $\bm{S_1} \bm{S_2} = \bm{U} \circ (\bm{S_1} \circ \bm{S_2})$ where $U$ is a fixed class, one can convert $E$ into a commutative Moufang loop as defined in \ref{CML} (for brevity, we shall speak of $CML\; E$). 
For all preceding examples of cubic surfaces where admissible equivalence was computed (such as those in Section \ref{admexamples} below) it led to associative commutative Moufang loops, i.e. to Abelian groups. 
The main result of this paper is constructing a cubic surface $V$ with non-associative $CML\; E$ for some admissible equivalence $\mathcal{A}$ on $V(k)$. 
This answers a long standing problem that was formulated after Theorem in Section 4.12 in \cite{manin1} in different but equivalent terms of existence of non-abelian quasi-groups  associated with some cubic surface over some field. In this paper the question about existence of cubic surfaces for which there exist non-associative structures is considered following \cite{manin3}   (see this formulation in Problem 11.11 in the Chapter II).
In our example (Section \ref{sectDef}) a cubic surface with non-associative Moufang loop of point classes is defined by
\begin{equation}
\label{diag}
    T_0^3+T_1^3+T_2^3+\theta T_3^3=0
\end{equation}
where $T_0, T_1, T_2, T_3$ is a homogenuous coordinate system in a projective space
$\mathbb{P}^3_k$ and $k=\mathbb{Q}_3(\theta)$, $\theta^3=1$ (quadratic extension of 3-adic numbers). The non-associativity Theorem \ref{maintheorem} is proved in Sections \ref{compsection}, \ref{proofsection}.\\
It was also can be shown (the work by the author in preparation) that in this example the admissible equivalence is universal, i.e. it is the finest possible admissible  equivalence on this surface.
Constructing an example of a cubic surface with a non-associative Moufang loop of point classes is not a trivial task, since one can show that for a cubic surface over a local field with a residue field of characteristic not equal 3 admissible equivalence on this surface over the local field generates an associative commutative Moufang loop of point classes (\cite{dimitrikanevsky3}).
And if a cubic surface defined over a local field with a residue field of characteristic 3 and certain non-degenerated conditions for V hold then a Moufang loop of point classes for this surface consists of only one element (Introduction, Theorem 2, \cite{swinnertondyer}). This leaves only a cubic surface over a local field with a residue field of characteristic 3 with "bad" reduction modulo a uniformizing element of the local field as a source of non-associative Moufang loops of point classes.

Our work opens a new direction in geometry of cubics with non-associative Moufang loop of point classes.
We can examine the asymptotic behavior of a specific subset of rational points on cubic surfaces over number fields, a subset generated by four rational points through the process of taking secants and tangents. A potential question to explore could be: Does this asymptotic behavior depend on whether this subset generates a nonassociative Moufang loops of point classes?
See a discussion in Section \ref{conjecture} on what classes of cubic surfaces over number fields may have non-associative Moufang loops of point classes. One can try to classify Moufang loops associated with these surfaces.

We would like to thank Prof. Vyacheslav Shokurov for reading the manuscript and suggesting improvements to it.

\section{CH-quasigroups, Moufang loops}
\numberwithin{equation}{section}
\label{summary}
Here, we give a brief survey of some algebraic and geometric structures which emerge in the theory of cubic hypersurfaces, initiated by Yu. Manin  (\cite{manin3}), Chapters I and II).
\begin{defi}
\label{quasi}
    A set $E$ with a binary composition law $E \times E \rightarrow E: (x,y) \mapsto x \circ y$ is called a symmetric quasigroup if it satisfies one of the following equivalent conditions:\\
    (i) The three-place relation $L(x,y,z): x \circ y = z$ is invariant under all permutations of $x,y,z$\\
    (ii) The following identities hold for all $x,y \in E:$
    $$x \circ y = y \circ x $$
    $$x \circ (x \circ y) = y$$
\end{defi}
A set $E$ of non-singular points on an irreducible cubic curve $V$ embedded in a projective plane over a field $k$ with a collinearity relationship $\{(x,y,x\circ y)\} = L \subset E \times E \times E$ provides an example of a symmetric quasigroup in the Definition \ref{quasi}.
One can choose a fixed element $u \in E$ and introduce on $E$ the new composition law: $xy = u \circ (x \circ y)$ for $x,y \in E$. Then $E$ becomes an Abelian group with $u$ as unit element. Generalization of this example leads to the following:
\begin{defi}
\label{abel}
    A symmetric quasigroup $E$ is called Abelian if it satisfies one of the following equivalent conditions:\\
    (i)  There exists an Abelian group structure on $E$ with composition law $(x,y) \mapsto xy$, and there is an element $c \in E$ such that $x \circ y = c x^{-1} y ^{-1}$ for all $x, y \in E$.\\
    (ii) For any element $u \in E$ the composition law $xy = u \circ (x \circ y)$ turns $E$ into an Abelian group.
\end{defi}
Let us now assume that $V$ is an irreducible cubic hypersurface of dimension $ \geq 2$ embedded in some projective space; let $E$ be the set of non-singular points on $V$. Let $L(x,y,z)$ be the three-place relationship on $E$ defined by the condition "$x,y,z$ are collinear". It is symmetric, but since in general it does not define a binary composition we consider a quotient set of $E$ such that induced relation of "collinearity" comes from a symetric quasigroup composition law. Since any three points of $V$ are contained in the intersection of $V$ with a plane, that defines a cubic curve we obtain the result that any set of three elements of the quasigroup generates an Abelian quasigroup.
\begin{defi}
\label{CH}
    A CH-quasigroup (CH stands for Cubic Hypersurface) is a symmetric quasigroup in which any three elements generate an Abelian subquasigroup.
\end{defi}
Next, by analogy with Abelian case we have the following
\begin{defi}
\label{CML}
    A set $E$ with composition law $(x,y) \mapsto xy = u \circ (x \circ y)$, where $u$ is some fixed element, is called a commutative Moufang loop (CML) if it is satisfied the following axioms:\\
    (i) Commutativity: $xy = yx$ for all $x,y \in E$.\\
    (ii) Unit element: $ux = x$ for all $x \in E$.\\
    (iii) Inverses: there exists a map $E \rightarrow E: x \mapsto x^{-1}$ such that $x^{-1}(xy) = y$ for all $x,y \in E$\\
    (iv) Weak associativity:\\
        
        for three factors: 
        $$x(xy) = x^2y;$$
        
        for four factors: 
        $$(xy)(xz) = x^2(yz),$$
        $$x(y(xz)) = (x^2y)z. $$
\end{defi}
The connection between CH-quasigroups and CMLs is given in the following (\cite{manin3}, Theorem 1.5)
\begin{theorem}
If $E$ is a CH-quasigroup and $u \in E$, then the composition law $xy = u\circ(x \circ y)$ turns $E$ into a CML.
\end{theorem}
Next, we describe some general structural results  for sets of classes of points on cubic hypersurfaces obtained by Yu. Manin.
\begin{defi}
\label{gen}
    Let $V \subset P^n$ be a cubic hypersurface. A geometric point $x$ on $V$ is called a point of general type if the following conditions are satisfied:\\
    (i) $x$ is non-singular point on $V$.\\
    (ii) The intersection $C(x)$ of the tangent hyperplane at $x$ with $V$ is geometrically irreducible and reduced.\\
    (iii) The singular point $x \in C(x)$ is not conical on $C(x)$. 
\end{defi}
\begin{theorem}
\label{period6}
Let $V$ be a cubic hypersurface of dimension $\geq 2$ having a point of general type, let $V(k)$ be the set of non-singular $k$-points. Let $\mathcal{A}$ be some admissible equivalence relation on $V(k)$ and $E = V(k)/\mathcal{A}$ be the corresponding symmetric quasigroup. Then the following statements hold:\\
(i) Every equivalence class is dense in the Zariski topology. If $k$ is algebraically closed, then all points of $V(k)$ are equivalent to one another.\\
(ii) $E$ is CH-quasigroup.\\
(iii) Let introduce on $E$ CML by means of the composition law $XY = U \circ (X \circ Y)$, where $U \in E$ is some fixed class. Then the relation $X^6 = 1$ holds on $E$.
\end{theorem}
From general structural results described in Chapter 1 in \cite{manin3} it also follows that under conditions of Theorem \ref{period6} $E$ is the direct product of an Abelian group of exponent 2 and a CML of exponent 3; and that if $E$ has a finite number of generators, then it is finite and consists of $2^a3^b$ elements.
\section{Some prior results on admissible equivalences}
\label{admexamples}
\subsection{R-equivalence} 
\label{requiv}
$R$-equivalence on points $V(k)$ of a cubic hypersurface $V$ over a filed $k$ is defined as the weakest equivalence relation such that for any two points $P, Q$ the following holds. There exists a sequence of points $T_0, T_1, ... T_n \in V(k)$ such that $T_0=P, T_n=Q$ and points $T_i$ and $T_{i+1}$ are directly $R$-connected, i.e. there exist a k-morphism $f: \mathbb{P}_k^1 \rightarrow V$  such that $f(0) = T_i$ and $f(\infty)=T_{i+1}$.
The $R-equivalence$ can be computed efficiently for Chatelet surfaces  for which it leads to associative commutative Moufang loops of point classes (\cite{manin3}, Proposition 45.10).
\subsection{Cubic surfaces over finite and local fields}
\label{lockan}
Let $k$ be a $\mathfrak{p}$-adic field and $\tilde k$ is residue field. We denote reduction $\mod \mathfrak{p}$ by a tilde. 
We take $V$ to be nonsingular with equations $F(X, Y, T, Z)=0$, where the coefficients of $F$ are $\mathfrak{p}$-adic integers; 
then $\tilde V$ is given by $\tilde F=0$. Theorem 2 in \cite{swinnertondyer} states that if $\tilde V$ smooth and Char $\tilde k \neq 2,3$ 
or some non degeneration conditions hold in Char $\tilde k =2,3$  then all the rational points of $V$ belong to the same class. 
Some cases for Char $\tilde k=2$ were considered in \cite{dimitrikanevsky} for which  it was shown that CML associated with these cases 
is associative. Our Theorem \ref{maintheorem} in this paper 
falls into the case Char $\tilde k=3$ and $\tilde V$ is non-smooth.
\subsection{Brauer equivalence} 
\numberwithin{equation}{section}
\label{brauer}
(\cite{manin3}, Chapter VI).
The Brauer equivalence on a set of rational points on an algebraic variety $V$ over a filed $k$ is defined via pairing $V(k) \times Br(V) \rightarrow Br(k)$, where $B=Br(V)$ is the Brauer group of $V$ and $Br(k)$ is the Brauer group of $k$. Namely, two points $P, Q \in V(k)$ are defined as $B$ equivalent if for any $b \in B$ we have for their specializations $b(P)=b(Q)$. The Brauer equivalence is admissible for cubic hypersurfaces $V$ and $CML\; V(k)/B$ is associative. 
\section{Definitions and main result}
\numberwithin{equation}{section}
\label{sectDef}
Let $K = \mathbb{Z}_3[\theta] \subset k$ be a quadratic extension of  a ring of 3-adic integers $\mathbb{Z}_3$ and let $\mathfrak{p} = 1 - \theta$. Then $3=-\theta^2 \mathfrak{p}^2$, $\mathfrak{p}$ generates the maximal ideal in $K$ such that factor of $K$ modulo this ideal is a finite field consisting of three elements.
We shall denote the additive valuation on $k$ by $\nu$, so that for $x \neq 0 \in k$ the exact power of $\mathfrak{p}$ divides $x$ is $\mathfrak{p}^{\nu(x)}$.\\
Let $\bm{P}=(t_0,t_1,t_2,t_3) \in V(k)$. We define $\bm{P} \mod \mathfrak{p}^n$ as $(t_0', ...t_3')$ where $t_i' = x t_i \mod \mathfrak{p}^n$, for some $x \in k$, $\nu(x t_i) \geq 0$ for all $i=0,...3$, and there is such $i \in\{0,1,2,3\}$ that
the following holds:
\begin{equation*}
    \label{canon}
    xt_i=1 \in K \text{   and for all  } j < i\;\; \nu(xt_j) >0
\end{equation*}
\begin{theorem}
\label{maintheorem}
Let $V$ be a cubic surface defined over $k$ by the equation (\ref{diag}). Let us define for $\bm{P},\bm{P'} \in V(k)$  $\bm{P} \sim \bm{P'} \mod \mathcal{A}$ if $\bm{P} = \bm{P'} \mod \mathfrak{p}^3$. Then $\mathcal{A}$ is admissible,
$CML\; (V(k)/\mathcal{A})$ is non-associative and consists of 243 elements.
\end{theorem}
\section{A non-associative composition example $\mod \mathfrak{p}^3$}
\numberwithin{equation}{section}
\label{compsection}
Here we give examples of points on $V(k)$ (defined by the equation (\ref{diag})) whose composition $\mod \mathfrak{p}^3$ does not obey  the associative law (in a sense of composition in $CML\; (V(k)/\mathcal{A}))$  in Theorem (\ref{maintheorem}).
\begin{lemma}
\label{lem1}
The tuple $(1,-1+\mathfrak{p}^2,\mathfrak{p},-\mathfrak{p})$ can be lifted to a point on $V(k)$.
\end{lemma}
\begin{lemma}
\label{lem2}
Let us define the following points in V(k)
\begin{equation*}
\label{U0}
    \bm{U_0} = (1,-1,0,0)
\end{equation*}
\begin{equation*}
\label{U1}
    \bm{U_1} = (1,0,-1,0)
\end{equation*}
\begin{equation*}
\label{U2}
    \bm{U_2} = (0, 1,-1,0)
\end{equation*}
Then for any $\bm{R}=(x,y,z,u) \in V(k)$ the following triples of points are collinear:
\begin{equation}
\label{frst}
\{\bm{U_0}, \bm{R}, (y,x,z,u)\}
\end{equation}
\begin{equation}
\label{scnd}
\{\bm{U_1}, \bm{R}, (z,y,x,u)\}
\end{equation}
\begin{equation}
\label{thrd}
\{\bm{U_2}, \bm{R}, (x,z,y,u)\}
\end{equation}
\end{lemma}
\begin{proposition}
\label{prop1}
Let $\bm{Q_0} = (1,-1+\mathfrak{p}^2,\mathfrak{p},-\mathfrak{p}) + \mathfrak{p}^3 v$ be a point on $V(k)$ for some $v \in K^4$.\\
Let us denote the following points on $V(k)$ as:
\begin{align*}
    \bm{Q_1} = (1,0,-1,0)\\
    \bm{Q_2} = (0,1,-\theta,0)
\end{align*}
And let's for any collinear triple of points $\bm{R_0}, \bm{R_1}, \bm{R_2} \in V(k)$, $\bm{R_0} \neq \bm{R_1}$  define $\bm{R_0}+\bm{R_1}$ as a point on $V(k)$ such that
$(\bm{U_0}, \bm{R_2}, \bm{R_0}+\bm{R_1})$ are collinear (this definition makes sense since there are no straight lines on $V$ defined over $k$). Then the following non-associative inequality holds:
\begin{equation*}
    (\bm{Q_0}+\bm{Q_1})+\bm{Q_2} \neq \bm{Q_0}+(\bm{Q_1}+\bm{Q_2}) \mod \mathfrak{p}^3
\end{equation*}
\end{proposition}
Since Theorem \ref{maintheorem} implies that the equivalence relation $\mod \mathcal{A}=\mod \mathfrak{p}^3$ on $V(k)$ is admissible Proposition \ref{prop1} gives examples of a triple of classes in $CML\; (V(k)/\mathcal{A}))$ with non-associative composition law.
\subsection{Remark}
\numberwithin{equation}{section}
The points $\bm{U_i}$ in Lemma \ref{lem2} are Eckhardt points, since if considered over the algebraic closure of $k$ a tangent plane to any of these points intersects $V$ in a curve that decomposes into three straight lines passing through this point. It is well known that any Eckhardt point on any projective cubic surface induces an automorhism of this surface that is a restriction of a reflection in a 3 dimensional projective space that contains this surface.  
\subsection{Proof of Lemma \ref{lem1}}
\numberwithin{equation}{section}
Let's show that for some $y \in K$:
\begin{equation}
\label{Q}
    (1,-1+\mathfrak{p}^2+\mathfrak{p}^3y,\mathfrak{p},-\mathfrak{p}) \in V(k)
\end{equation}
Substituting (\ref{Q}) into (\ref{diag})  we get:
\begin{equation}
\label{this}
    3\mathfrak{p}^2-3\mathfrak{p}^4+\mathfrak{p}^6 +3(-1+\mathfrak{p}^2)^2\mathfrak{p}^3y+3(-1+\mathfrak{p}^2)\mathfrak{p}^6z^2+\mathfrak{p}^{9}z^3+\mathfrak{p}^4=0
\end{equation}
Substituting  $3=-\theta^2\mathfrak{p}^2$ into (\ref{this}), using $\mathfrak{p}^3-\theta \mathfrak{p}^3=\mathfrak{p}^4$ and dividing the equation by $\mathfrak{p}^5$ we get
\begin{equation*}
    (1+\theta) - \theta^2y + \mathfrak{p}( ... \text{ terms with coefficients in K }...) = 0
\end{equation*}
The result now follows from Hensel's lemma.\\
\subsection{Proof of Lemma \ref{lem2}}
\numberwithin{equation}{section}
Without loss of generality we can assume that $\bm{R}=(1,y,z,u)$ and $\bm{U}_0 \neq \bm{R}$.
Taking a straight line through $\bm{U_0}, \bm{R}$ we get 
\begin{equation}
\label{u0r}
\bm{U_0}\circ \bm{R} = (1,-1 + (y+1)\tau, z\tau, u\tau)
\end{equation}
for some $\tau \in k$. Substituting (\ref{u0r}) into (\ref{diag}) and solving this equation in $\tau$ gives us:
\begin{align*}
3(y+1)-3(y+1)^2\tau +(y+1)^3\tau^2+z^3\tau^2+\theta u^3\tau^2=0\\
3(y+1)-3(y+1)^2\tau +(3y+3y^2)\tau^2=0
\end{align*}
This implies:
$$\tau = \frac{3(y+1)}{3y+3y^2} =\frac{1}{y}$$
and 
$$\bm{U_0}\circ \bm{R} =
\left (1,-1+(y+1)\frac{1}{y},z\frac{1}{y},u\frac{1}{y}\right ) = (y, 1, z, u) \in V(k)$$ 
Q.E.D.
Similarly one can prove (\ref{scnd}) and (\ref{thrd}).
\subsection{Proof of Proposition \ref{prop1}}
\numberwithin{equation}{section}
Since $(\bm{Q_0}+\bm{Q_1})+\bm{Q_2} \neq \bm{Q_0}+(\bm{Q_1}+\bm{Q_2})$ iff $(\bm{Q_0}+\bm{Q_1})\circ \bm{Q_2} \neq \bm{Q_0}\circ (\bm{Q_1}+\bm{Q_2}) $
we compute $(\bm{Q_0}+\bm{Q_1})\circ \bm{Q_2}$  and $\bm{Q_0}\circ(\bm{Q_1}+\bm{Q_2})$
\subsubsection{Additional notations and definitions}
\numberwithin{equation}{section}
Throughout this paper we will use the following notations:
\begin{itemize}
    \item For any tuples $v, v' \in K^n$ we set $v\equiv_i v'$ if $v = v' \mod \mathfrak{p^i}$
    \item $\tilde K_i :=K \mod \mathfrak{p^i}$
    \item For any tuples $v, v' \in \tilde K^n_i$ we set
$v \sim v' $ if there is $x \in K,\; \nu(x) =0$ such that $v = xv'$ (with the obvious multiplication action of $x$ on $\tilde K_i^n$).
\item For any tuples $v=(t_0,t_1,t_2, t_3), v'=(t'_0,t'_1,t'_2,t'_3) \in K^4$ that define the same point in  $V(K)$ we write $v=v'$.
\end{itemize}
\subsubsection{\texorpdfstring{Computation of $(\bm{Q_0}+\bm{Q_1})\circ \bm{Q_2}$}%
{Computation of (Q_0+Q_1)\circ Q_2}}
$\bm{Q_1} = \bm{U_1}$ from Lemma {\ref{lem2}} interchanges a first and third coordinates when it is applied to $\bm{Q_0}$ and therefore
\begin{equation*}
    \bm{Q_0} \circ \bm{Q_1} \mod \mathfrak{p}^3 \sim (\mathfrak{p},-1+\mathfrak{p}^2,1,-\mathfrak{p})
\end{equation*}
Next, $\bm{U_0}$ interchanges a first and second coordinates in $\bm{Q_0} \circ \bm{Q_1}$ and therefore 
\begin{equation*}
    \bm{U_0}\circ (\bm{Q_0}\circ \bm{Q_1})=\bm{Q_0} + \bm{Q_1} \equiv_3 (-1+\mathfrak{p}^2,\mathfrak{p},1,-\mathfrak{p})
\end{equation*}
Next, we compute $(\bm{Q_0} + \bm{Q_1}) \circ \bm{Q_2}$.
We have $\bm{Q_2} \mod \mathfrak{p}^3 = (0,1,-\theta,0) \sim (0,-\theta^2,1,0)$\\
Therefore,
\begin{equation}
\label{tau}
(\bm{Q_0} + \bm{Q_1}) \circ \bm{Q_2} = ((-1+\mathfrak{p}^2)\tau , -\theta^2 + (\mathfrak{p}+\theta^2) \tau , 1, - \mathfrak{p}\tau) + \mathfrak{p}^3v' \tau
\end{equation}
for some $\tau \in k$, $v' \in K^4$.
Substituting (\ref{tau}) in (\ref{diag}) and solving a quadratic equation in $\tau$ we get
\begin{equation}
\label{tau1}
    \tau = \frac{\theta^2+x\mathfrak{p}^3}{\mathfrak{p} + y\mathfrak{p}^3 }
\end{equation}
for some $x,y \in K$.
Then multiplying all terms in (\ref{tau}) by $\mathfrak{p}$, substituting (\ref{tau1}) into (\ref{tau}), 
and taking modulo $\mathfrak{p}^3$ 
we get
\begin{align*}
(\bm{Q_0} + \bm{Q_1}) \circ \bm{Q_2} \mod \mathfrak{p}^3 =\\ (-\theta^2+\mathfrak{p}^2, \theta , \mathfrak{p} , -\mathfrak{p} \theta^2) \sim (1,-\theta^2 - \mathfrak{p}^2 , - \mathfrak{p} + \mathfrak{p}^2 , \mathfrak{p}) \in \tilde K^4_3
\end{align*}
\subsubsection{\texorpdfstring{Computation of $\bm{Q_0}\circ (\bm{Q_1}+\bm{Q_2})$}%
{Computation of Q_0\circ (Q_1+Q_2)}}
We have 
$$\bm{Q_1} \circ \bm{Q_2} = (-\theta,1,0,0) $$
$$\bm{Q_1} + \bm{Q_2} = (1:-\theta:0:0) $$
\begin{equation*}
\bm{Q_0} \circ (\bm{Q_1} + \bm{Q_2}) = 
\end{equation*}
\begin{equation*}
    ((1,-1+\mathfrak{p}^2,\mathfrak{p},-\mathfrak{p})+\mathfrak{p}^3v) \circ (1,-\theta,0,0)=(1, -\theta + (-1 + \mathfrak{p}^2 + \theta)\tau , \mathfrak{p}\tau , -\mathfrak{p}\tau) \mod \mathfrak{p}^3 
\end{equation*}    
$$\sim (1, -\theta^2 - \mathfrak{p}^2, -\mathfrak{p} + \mathfrak{p}^2, \mathfrak{p} - \mathfrak{p}^2) $$
where  $\tau = -\frac{1}{1+\mathfrak{p}} \mod \mathfrak{p}^3$.\\
Proposition \ref{prop1} now follows from the fact that
$$ (1,-\theta^2 - \mathfrak{p}^2 , - \mathfrak{p} + \mathfrak{p}^2 , \mathfrak{p}) \neq (1, -\theta^2 - \mathfrak{p}^2, -\mathfrak{p} + \mathfrak{p}^2, \mathfrak{p} - \mathfrak{p}^2) \mod \mathfrak{p}^3$$
\subsection{Number of points}
\numberwithin{equation}{section}
In this subsection we compute number of points on $V(K)$ modulo $\mathfrak{p^i}$, $i=1,2,3$ and also give some representation for elements in these sets that will be needed for calculations in  the Section \ref{proofsection}.
\subsubsection{Number of points modulo $\mathfrak{p}$}
\numberwithin{equation}{section}
\begin{lemma}
$V(K) \mod \mathfrak{p}$ consists of three tuples in $\mathbb{F}_3^4$, where $\mathbb{F}_3$ is a field of 3 elements.
\begin{align*}
(1,-1,0,0)\\
(1,0,-1,0)\\
(0,1,-1,0)
\end{align*}
\end{lemma}
\textbf{Proof}\\
First, we have that after permutation of coordinates any canonical representation of a solution of (\ref{diag}) has the following form:
\begin{equation}
\label{solap}
(X,Y,Z,U)=(1, - 1 + \mathfrak{p} \tilde{y}, \mathfrak{p} \tilde{z}, \mathfrak{p} \tilde{u})
\end{equation}
where $ \tilde{y}, \tilde{z},  \tilde{u} \in K$.\\
Indeed, 
let $X=1$, $Y= a + \mathfrak{p}b$, $Z = c + \mathfrak{p}d$, $U=e + \mathfrak{p}f$ where $a,c,e \in \{-1,0, 1\}$ and $b,d,f \in K$.\\
Putting this into (\ref{diag}) gives:
\begin{align*}
\label{modpp}
1 + (a + \mathfrak{p}b)^3 + (c + \mathfrak{p}d)^3 + \theta(e + \mathfrak{p}f)^3) = \\
1 + a^3 + 3\mathfrak{p}a^2b + 3\mathfrak{p}^2ab^2 + \mathfrak{p}^3b^3 + c^3 + 3\mathfrak{p}c^2 d+ 3\mathfrak{p}^2cd^2 + \mathfrak{p}^3d^3 + \theta e^3 + 3\mathfrak{p}\theta e^2f + 3\mathfrak{p}^2\theta ef^2+\mathfrak{p}^3\theta f^3 = 0
\end{align*}
This gives
\begin{equation}
\label{modap}
1 + a^3 + c^3 + \theta e^3 = 0 \mod  \mathfrak{p}^3
\end{equation}
This implies $e = 0 $ in (\ref{modap})) (i.e. $U = 0 \mod  \mathfrak{p}$ in (\ref{solap})) and 
\begin{equation*}
\label{modap1}
1 + a^3 + c^3 = 0 \mod  \mathfrak{p}^3
\end{equation*}
Substituting all choices of $a,c \in \{-1,0, 1\}$ gives two solutions:
$a = -1, \; c = 0$ or $a = 0, \; c = -1$ Q.E.D.
\subsubsection{Number of points modulo $\mathfrak{p}^2$}
\numberwithin{equation}{section}
\begin{lemma}
$V(K) \mod \mathfrak{p}^2$ consists of 27 tuples over $K \mod \mathfrak{p}^2$
\begin{align*}
(1,-\theta^i,\epsilon \mathfrak{p},-\epsilon \mathfrak{p})\\
(1,\delta \mathfrak{p},-\theta^j,-\delta \mathfrak{p})\\
(\sigma \mathfrak{p},1,-\theta^r,-\sigma \mathfrak{p})
\end{align*}
where $i,j,r,\epsilon, \delta, \sigma \in \{0,1,-1\}$
\end{lemma}
\textbf{Proof}\\
Let 
\begin{align*}
\tilde y = -1 + \mathfrak{p} y\\
\tilde z = \mathfrak{p} z\\
\tilde u = \mathfrak{p} u
\end{align*}
where $ z, u \in \{0,1,-1\} \subset K$ and $y \in K$.
Substituting $\tilde y,\tilde z,\tilde u$ into a cubic equation
$$1 +  \tilde y^3 +  \tilde z^3 + \theta  \tilde u^3 = 0 $$
gives us:
\begin{equation*}
3\mathfrak{p} y - 3\mathfrak{p}^2 y^2 + \mathfrak{p}^3 y^3 + \mathfrak{p}^3 z^3 + \theta \mathfrak{p}^3 u^3 = 0 
\end{equation*} 
or
\begin{equation*}
y - \mathfrak{p} y^2  - \theta y^3 -\theta z^3 - \theta^2 u^3 = 0 
\end{equation*} 
Since $y - \mathfrak{p} y^2  - \theta y^3 = 0 \mod \mathfrak{p}$ we have $ z +  u = 0 \mod \mathfrak{p}$ or 
$ z = - u$. Q.E.D. 
\subsubsection{Number of points modulo $\mathfrak{p}^3$}
\numberwithin{equation}{section}
Let's define polynomials maps: $\{\mathbb{Z},K^4\} \rightarrow K^4$:
\begin{align}
\label{param}
P(i,\epsilon, y,z,u) = (1,-\theta^i + \mathfrak{p}^2y,\epsilon \mathfrak{p} + \mathfrak{p}^2z,-\epsilon \mathfrak{p} + \mathfrak{p}^2 u)\\
\label{param1}
Q(j,\delta, y_1, z_1, u_1)=(1,\delta \mathfrak{p} +  \mathfrak{p}^2 y_1,-\theta^j + \mathfrak{p}^2 z_1,-\delta \mathfrak{p} + \mathfrak{p}^2 u_1)\\
\label{param2}
R(r, \sigma, x, z_2, u_2)=(\sigma \mathfrak{p} + \mathfrak{p}^2 x,1,-\theta^r + \mathfrak{p}^2 z_2,-\sigma \mathfrak{p} + \mathfrak{p}^2 u_2)
\end{align}
where $i,j,r$ are integers and $\epsilon, \delta, \sigma,  z, z_1, z_2, u, u_1, u_2, x, y, y_1 \in K$.
And lets define subsets:
\begin{align*}
\Lambda_P=\{P(i,y, y,z,u)\}\\
\Lambda_Q= \{Q(j,z_1, y_1, z_1, u_1)\}\\
\Lambda_R= \{R(r, z_2, x, z_2, u_2)\}\\
\Lambda = \Lambda_P\cup\Lambda_Q\cup\Lambda_R
\end{align*}
where $i,j,r,\epsilon, \delta, \sigma,  z, z_1, z_2, u, u_1, u_2, x, y, y_1 \in  \{0,1,-1\} \subset K$.
\begin{lemma}
\label{243}
$V(K) \mod \mathfrak{p}^3 \cong \Lambda $ and $\#\Lambda = 243$.
\end{lemma}
\textbf{Proof}\\
Substituting $(1, \tilde y, \tilde z, \tilde u) := P(i,\epsilon,y',z,u)$, where $y' \in K$  into a cubic equation with affine coordinates $\tilde y, \tilde z, \tilde u$ 
$$1 + \tilde y^3 + \tilde z^3 + \theta \tilde u^3 = 0$$
gives us
\begin{align*}
\epsilon \mathfrak{p}^4 + 3\theta^{2i} \mathfrak{p}^2y' - 3 \theta^i \mathfrak{p}^4 y'^2 + \mathfrak{p}^6 y'^3 +\\ 
3\epsilon^2 \mathfrak{p}^4 z + 3\epsilon \mathfrak{p}^5 z^2 + \mathfrak{p}^6 z^3 +\\
3\theta \epsilon^2 \mathfrak{p}^4 u - 3\theta \epsilon \mathfrak{p}^5 u^2 + \theta \mathfrak{p}^6 u^3 = 0 
\end{align*}
and
\begin{align*}
\epsilon  - \theta^{2i+2} y'  - \mathfrak{p}^2 \theta^2 \epsilon^2  z -  \mathfrak{p}^2 \epsilon^2 u + \\ 
 \mathfrak{p}^2 \theta^{2+i} y'^2 + \mathfrak{p}^2 y'^3 +\\ 
 - \mathfrak{p}^3 \theta^2\epsilon z^2 + \mathfrak{p}^2 z^3 +\\
\mathfrak{p}^3  \epsilon  u^2 + \theta \mathfrak{p}^2 u^3 = 0
\end{align*}
For any $z', u' \in K$  there is a solution for $y' \in K$ of this equation according to Hensel's lemma 
\begin{align*}
\epsilon  - \theta^{2i+2} y'  - \mathfrak{p}^2 \theta^2 \epsilon^2  z' -  \mathfrak{p}^2 \epsilon^2 u'\\ 
 \mathfrak{p}^2 \theta^{2+i} {y'}^2 + \mathfrak{p}^2 {y'}^3 +\\ 
 - \mathfrak{p}^3 \theta^2\epsilon {z'}^2 + \mathfrak{p}^2 {z'}^3 +\\
\mathfrak{p}^3  \epsilon  {u'}^2 + \theta \mathfrak{p}^2 {u'}^3 = 0
\end{align*}
This equation implies 
\begin{equation*}
y' = \theta^{i+1} \epsilon + \text{ terms in K divisible by  }  \mathfrak{p}^2
\end{equation*}
This implies for $y'=y\mod \mathfrak{p}$:
\begin{equation*}
  \epsilon = y  
\end{equation*}
There are 9 choices for $\{z,u\}$ such that $z = z' \mod \mathfrak{p}, u = u' \mod \mathfrak{p}$ in the above equation. Each of them uniquely defines $\epsilon - \theta^{2i+2}y' = \epsilon - y = 0 \mod \mathfrak{p}$ where $y=y' \mod \mathfrak{p}$. There are three choices for $\epsilon, y$ such that $\epsilon - y = 0$. There are 3 choices for $i$. \\
I.e.  $\#\Lambda_P = 81$.\\
Repeating similar calculations for $\Lambda_Q$ and $\Lambda_R$ gives us
\begin{equation*}
  \delta = z_1  
\end{equation*}
\begin{equation*}
  \sigma = z_2  
\end{equation*}
And finally we get: $\#\Lambda = 3*81= 243$. 
Q.E.D.
\subsubsection{Remark}
\numberwithin{equation}{section}
Using parametric representation for elements in $\Lambda$ like in (\ref{param}, \ref{param1}, \ref{param2}) and for composition of points (\ref{casea}) one can give a full description of triples of points on $V(k)$ that leads to non-associative composition relations and describe the associative center in $CML\; (V(k)/\mathcal{A}))$ (the work by the author in preparation).
\section{Proof of admissibility of relation modulo $\mathfrak{p}^3$ on $V(k)$}
\numberwithin{equation}{section}
\numberwithin{equation}{section}
\label{proofsection}
In the notation of Theorem \ref{maintheorem} we prove:
\begin{proposition}
\label{adm}
If $\bm{P_1}, \bm{P_2}, \bm{P_3} \in V(k)$ and $\bm{P_1'}, \bm{P_2'}, \bm{P_3'} \in V(k)$ are collinear, $\bm{P_1} \equiv_3 \bm{P_1'}$,
$\bm{P_2} \equiv_3 \bm{P_2'}$ then $\bm{P_3} \equiv_3 \bm{P_3'}$.
\end{proposition}
We prove some statements in the following three cases and use them to prove Proposition \ref{adm}:\\
    \textbf{Case 1}:  $\bm{P_1} = \bm{P_2}$, $\bm{P_1'} = \bm{P_2'}$ (points on a tangent plane).\\
    \textbf{Case 2}:  $\bm{P_1} \neq \bm{P_2} \mod \mathfrak{p}$\\
    \textbf{Case 3}: $\bm{P_1} = \bm{P_2} \mod \mathfrak{p}$
\subsection{Case 1: Points on a tangent plane}
\begin{proposition}
\label{propap}
Let ${\bm P} , {\bm P'} \in V(k)$ and ${\bm P'}$  is lying on a tangent plane to $V$ at ${\bm P}$, then ${\bm P} \equiv_3 {\bm P'} $
\end{proposition}
\textbf{Proof}\\
Let ${\bm P} = (1, - 1 + \mathfrak{p} \tilde{y}, \mathfrak{p} \tilde{z}, \mathfrak{p} \tilde{u})$ . Substituting this into  (\ref{diag}) gives 
\begin{align*}
3\mathfrak{p}\tilde{y} - 3 \mathfrak{\mathfrak{p}}^2 \tilde{y}^2 + \mathfrak{p}^3\tilde{y}^3 + \mathfrak{p}^3 \tilde{z}^3 + \theta \mathfrak{p}^3 \tilde{u}^3 = 0 \\
-\theta^2\mathfrak{p}^3\tilde{y} + \theta^2 \mathfrak{p}^4 \tilde{y}^2 + \mathfrak{p}^3 \tilde{y}^3 + \mathfrak{p}^3 \tilde{z}^3 + \theta \mathfrak{p}^3 \tilde{u}^3 = 0\\
\end{align*}
or 
\begin{equation}
\label{eqsol}
-\theta^2\tilde{y} + \theta^2 \mathfrak{p} \tilde{y}^2 +  \tilde{y}^3 +  \tilde{z}^3 + \theta  \tilde{u}^3 = 0
\end{equation}
We can interpret  (\ref{eqsol}) as solution of the following equation in variables $\tilde Y, \tilde Z, \tilde U$:
\begin{equation}
\label{eq}
-\theta^2{\tilde Y} + \theta^2 \mathfrak{p} {\tilde Y}^2 +  {\tilde Y}^3 +  {\tilde Z}^3 + \theta  {\tilde U}^3 = 0
\end{equation}
Let us consider a transformation:
\begin{equation*}
 \mathcal{T}: (1,Y,Z,U) \mapsto \left (\frac{Y+1}{\mathfrak{p}}, \frac{Z}{\mathfrak{p}}, \frac{U}{\mathfrak{p}}\right)
\end{equation*}
Then
\begin{equation}
\label{pointap1}
\mathcal{T}{\bm P} = \left ( \tilde{y} ,  \tilde{z}, \tilde{u} \right)
\end{equation}
and
\begin{equation}
\label{pointap2}
\mathcal{T}{\bm P'} = \left (\tilde{y} + \tilde{\alpha}t,  \tilde{z} + \tilde{\beta}t, \tilde{u} + \tilde{\gamma} t \right)
\end{equation}
for some $( \tilde{\alpha},  \tilde{\beta},  \tilde{\gamma}) \in K^3$ such that 
\begin{equation}
\label{ordap}
\nu(\tilde{\alpha})\nu(\tilde{\beta})\nu(\tilde{\gamma})=0
\end{equation} 
and $t \in k$.
Substituting (\ref{pointap2}) into (\ref{eq}) we get the following equalities:
\begin{equation*}
\label{eqap}
-\theta^2(\tilde{y} + \tilde{\alpha}t) + \theta^2 \mathfrak{p} (\tilde{y} + \tilde{\alpha}t)^2 +  (\tilde{y} + \tilde{\alpha}t)^3 +  (\tilde{z} + \tilde{\beta}t)^3 + \theta  (\tilde{u} + \tilde{\gamma} t)^3 = 0
\end{equation*}
or
\begin{equation*}
\label{compressap}
l_0 + l_1t + l_2t^2=0
\end{equation*}
where
\begin{equation*}
\label{eqap0}
l_0=-\theta^2\tilde{\alpha} +  \mathfrak{p}\theta^2 2 \tilde{y} \tilde{\alpha} + 3\tilde{y}^2\tilde{\alpha} + 3 \tilde{z}^2 \tilde{\beta} + 3\theta\tilde{u}^2 \tilde{\gamma}
\end{equation*}
\begin{equation*}
\label{eqap1}
l_1=\mathfrak{p}\theta^2\tilde{\alpha}^2 + 3 \tilde{y} \tilde{\alpha}^2  + 3 \tilde{z} \tilde{\beta}^2 + 3\theta\tilde{u} \tilde{\gamma}^2
\end{equation*}
\begin{equation*}
\label{eqap2}
l_2=\tilde{\alpha}^3 +  \tilde{\beta}^3 + \theta \tilde{\gamma}^3
\end{equation*}
Since $\bm P'$ belongs to an intersection of a tangent plane at $\bm P$ with $V$ we have $l_0=0$.
Therefore, $\nu(\tilde \alpha) \geq 2$ or $ \tilde \alpha = \mathfrak{p}^2 \bar \alpha$ for some $\bar \alpha \in K$.
Next, we have:
\begin{equation}
\label{tap}
t = -\frac{ \mathfrak{p}^5 \theta^2 \bar \alpha^2 +  3\mathfrak{p}^4 \tilde y  \bar \alpha ^2 + 3 \tilde z \tilde \beta^2 + 3 \theta \tilde u \tilde \gamma^2}{\mathfrak{p}^6\bar \alpha^3 + \tilde \beta^3 + \theta \tilde \gamma ^3}
\end{equation}
If $\tilde \beta^3 + \theta \tilde \gamma ^3 = 0 \mod \mathfrak{p}^2$  then $\nu(\tilde \beta) > 0$ and $\nu(\tilde \gamma) > 0$ that contradicts to a fact that $\nu(\tilde \alpha) > 0$ and assumption (\ref{ordap}).
If $\tilde \beta^3 + \theta \tilde \gamma ^3 \neq 0 \mod \mathfrak{p}$  then $\nu(t) \geq 2$ and therefore 
$\mathcal{T}{\bm P} \equiv_2 \mathcal{T}{\bm P'}$.
This implies that $\bm P \equiv_3 \bm P'$.
This leaves us the last case to consider:
\begin{equation*}
\label{ord1}
\nu(\tilde \beta^3 + \theta \tilde \gamma ^3) = 1
\end{equation*}
This gives as options for $ \tilde \beta = 1 \mod \mathfrak{p}, \; \tilde \gamma = -1 \mod \mathfrak{p}$ or $ \tilde \beta = -1 \mod \mathfrak{p}, \; \tilde \gamma = 1 \mod \mathfrak{p}$.
This allows to represent (\ref{tap}) as
\begin{equation}
\label{tap2}
t = -\frac{ \mathfrak{p}^5 \theta^2 \bar \alpha^2 +  3\mathfrak{p}^4 \tilde y  \bar \alpha ^2 + 3 (\tilde z  +   \tilde u - \mathfrak{p}\tilde u)  + p^3 g}{\mathfrak{p}^6\bar \alpha^3 + \tilde \beta^3 + \theta \tilde \gamma ^3}
\end{equation}
for some $g \in K$.
Note that in (\ref{eqsol}) 
\begin{equation*}
\label{part}
-\theta^2\tilde{y} + \theta^2 \mathfrak{p} \tilde{y}^2 +  \tilde{y}^3  = 0 \mod \mathfrak{p} 
\end{equation*}
Therefore 
\begin{equation}
\label{z3u3}
\tilde{z}^3 + \theta  \tilde{u}^3  = 0 \mod \mathfrak{p}
\end{equation}
Next:
\begin{align*}
\tilde{z}^3 + \theta  \tilde{u}^3  = \tilde{z}^3 +  \tilde{u}^3 -  \tilde{u}^3+\theta  \tilde{u}^3 = \tilde{z}^3 +  \tilde{u}^3 -  \mathfrak{p} \tilde{u}^3 = \\
 (\tilde{z} +  \tilde{u})^3 -  \mathfrak{p} \tilde{u}^3 - 3\tilde{z}^2\tilde{u} - 3\tilde{z}\tilde{u}^2
\end{align*}
This implies that 
\begin{equation*}
\label{sumzu3}
(\tilde{z} +  \tilde{u})^3 = 0 \mod p
\end{equation*}
or that 
\begin{equation}
\label{sumzu}
\tilde{z} +  \tilde{u} = 0 \mod p
\end{equation}
This also means that in (\ref{tap2}) 
\begin{align*}
\nu(\mathfrak{p}^5 \theta^2 \bar \alpha^2 +  3\mathfrak{p}^4 \tilde y  \bar \alpha ^2  + 3 (\tilde z  +   \tilde u - \mathfrak{p}\tilde u) + \mathfrak{p}^3 g) \geq 3 \\
\nu(\mathfrak{p}^6\bar \alpha^3 + \tilde \beta^3 + \theta \tilde \gamma ^3) = 1
\end{align*}
Therefore $\nu(t) \geq 2$, i.e. 
$\mathcal{T}{\bm P} \equiv_2 \mathcal{T}{\bm P'}$
and finally $\bm P \equiv_3 \bm P'$.
\subsection{\texorpdfstring{Case 2: $\bm{P_1} \neq \bm{P_2} \mod \mathfrak{p}$}%
{Case 2: P_1 \neq P_2 \mod \mathfrak{p}}}
Here we compute composition of points on $V(k)$ that lie above $P(i,\epsilon, y,z,u), Q(j,\delta, y_1, z_1, u_1)$. Namely, we show that
if $\bm{P}, \bm{Q}, \bm{R} \in V(k)$ are collinear and $\bm{P} =  P(i,\epsilon, y,z,u) \mod \mathfrak{p}^3$, $\bm{Q} = Q(j,\delta, y_1, z_1, u_1) \mod \mathfrak{p}^3$ then $\bm{R} = (R_0, R_1, R_2, R_3) \mod \mathfrak{p}^3 \in \tilde K_3^4$ where $R_0,R_1,R_2,R_3 $ are polynomials of $i,\epsilon, y,z,u, j,\delta, y_1, z_1, u_1$ with coefficients in $K$. This immediately gives us the Case 2 in the admissibility proof.

We set for $P(i,\epsilon,y,z,u)$:
\begin{align*}
    X=1\\
    Y = -\theta^i + \mathfrak{p}^2 y\\
    Z = \epsilon \mathfrak{p} + \mathfrak{p}^2 z = \mathfrak{p} y + \mathfrak{p}^2 z\\
    U = - \epsilon \mathfrak{p}+ \mathfrak{p}^2 u = - \mathfrak{p} y + \mathfrak{p}^2 u
\end{align*}
And we set for $Q(j,\delta, y_1, z_1, u_1) $
\begin{align*}
    X=1\\
    Y_1 = \delta \mathfrak{p} + \mathfrak{p}^2 y_1 = \mathfrak{p} z_1 + \mathfrak{p}^2 y_1\\
    Z_1 = -\theta^j + \mathfrak{p}^2 z_1\\
    U_1 = - \delta \mathfrak{p}+ \mathfrak{p}^2 u_1 = - \mathfrak{p} z_1 + \mathfrak{p}^2 u_1\\
\end{align*} 
And then:
\begin{align*}
    \alpha = Y_1-Y = \theta^i + \mathfrak{p} z_1 + \mathfrak{p}^2(y_1 - y)\\
    \beta = Z_1 - Z = -\theta^j - \mathfrak{p} y + \mathfrak{p}^2(z_1 - z)\\
    \gamma = U_1 - U = \mathfrak{p} (y - z_1) + \mathfrak{p}^2 (u_1 - u)
\end{align*}
Next we set:
\begin{align*}
    \tau = -\frac{l_0}{l_0+l_1} = -\frac{1}{1+l_1/l_0}\\
    l_0 = Y^2 \alpha + Z^2 \beta + \theta U^2 \gamma\\
    l_1 = Y \alpha^2 + Z \beta^2 + \theta U \gamma^2
\end{align*}
This gives us:
\begin{align*}
    \frac{l_1}{l_0}\equiv_3 \frac{\alpha Y^{-1} + Y^{-2}\alpha^{-1}Z\beta^2}{1+Y^{-2}\alpha^{-1}Z^2\beta} \equiv_3 \\
    (\alpha Y^{-1} + Y^{-2} \alpha ^{-1} Z \beta^2)(1 - Y^{-2} \alpha^{-1} Z^2 \beta) \equiv_3 \\
    \alpha Y^{-1} + Y^{-2} \alpha ^{-1} Z \beta^2 - Y^{-3} Z^2 \beta
\end{align*}
Next, we will compute components $\alpha Y^{-1}$, $Y^{-2} \alpha ^{-1} Z \beta^2 $ and $Y^{-3} Z^2 \beta$.
\begin{align*}
    Y^{-1} = \frac{1}{-\theta^i + \mathfrak{p}^2 y}= -\theta^{2i} \frac{1}{1 - \mathfrak{p}^2\theta^{2i}y} \equiv_3 \\
    -\theta^{2i}(1 + \mathfrak{p}^2y) \equiv_3 -\theta^{2i} - \mathfrak{p}^2y
\end{align*}
This gives us
\begin{align*}
    \alpha Y^{-1} = 
    [\theta^i + \mathfrak{p} z_1 + \mathfrak{p}^2(y_1 - y)](-\theta^{2i} - \mathfrak{p}^2y) \equiv_3\\
    -1 - \mathfrak{p}\theta^{2i}z_1 + \mathfrak{p}^2(y-y_1) - \mathfrak{p}^2y \equiv_3\\
    -1 - \mathfrak{p}\theta^{2i}z_1 - \mathfrak{p}^2 y_1
\end{align*}
Next we compute $Y^{-2} \alpha ^{-1} Z \beta^2 $
\begin{align*}
    Y^2 \equiv_3 \theta^{2i} + \mathfrak{p}^2y\\
    Y^{-2} = \frac{1}{\theta^{2i} + \mathfrak{p}^2y} \equiv_3 \theta^i \frac{1}{1 + \mathfrak{p}^2y} \equiv_3\\
    \theta^i(1 - \mathfrak{p}^2y) = \theta^i - \mathfrak{p}^2y\\
    \alpha^{-1} = \frac{1}{\theta^i + \mathfrak{p} z_1 + \mathfrak{p}^2(y_1-y) } \equiv_3\\
    \theta^{2i}\frac{1}{1 + \mathfrak{p}\theta^{2i}z_1 + \mathfrak{p}^2(y_1-y)} \equiv_3\\
    \theta^{2i} - \mathfrak{p}\theta^{i}z_1 - \mathfrak{p}^2(y_1-y) + \mathfrak{p}^2z_1^2
\end{align*}
This gives us:
\begin{align*}
    Y^{-2} \alpha ^{-1} Z \beta^2  \equiv_3 (\theta^i - \mathfrak{p}^2y)[\theta^{2i} - \mathfrak{p}\theta^{i}z_1 - p^2(y_1-y) + p^2z_1^2] (p y + p^2z) [-\theta^j - p y + p^2(z_1 - z)]^2 \equiv_3\\
    \theta^i(\theta^{2i} - \mathfrak{p} z_1) (\mathfrak{p} y + \mathfrak{p}^2z)(\theta^{2j} - \mathfrak{p} y) \equiv_3\\
    \theta^i(p\theta^{2i}y + \mathfrak{p}^2z - \mathfrak{p}^2y z_1)(\theta^{2j} - \mathfrak{p} y) \equiv_3\\
    \theta^i(p\theta^{2i+2j}y + \mathfrak{p}^2z - \mathfrak{p}^2y z_1 - \mathfrak{p}^2y^2) \equiv_3\\
     \mathfrak{p}\theta^{2j}y + \mathfrak{p}^2z - \mathfrak{p}^2y z_1 - \mathfrak{p}^2y^2 
\end{align*}
And next we get
\begin{equation*}
    Y^{-3} Z^2 \beta = \mathfrak{p}^2 y^2
\end{equation*}
Now we can compute:
\begin{align*}
    \frac{l_1}{l_0} \equiv_3 -1 - \mathfrak{p}\theta^{2i}z_1 - \mathfrak{p}^2 y_1 + \mathfrak{p}\theta^{2j}y + \mathfrak{p}^2z - \mathfrak{p}^2y z_1 - \mathfrak{p}^2y^2 - \mathfrak{p}^2 y^2 \equiv_3\\
    -1 + \mathfrak{p}(\theta^{2j}y - \theta^{2i} z_1) + \mathfrak{p}^2(-y_1 + z + y^2 - y z_1)
\end{align*}
This gives us finally:
\begin{equation*}
    \tau = -\frac{l_0}{l_0+l_1} = - \frac{1}{1+l_1/l_0} = -\frac{1}{\mathfrak{p}(\theta^{2j}y-\theta^{2i}z_1)+\mathfrak{p}^2(-y_1+z+y^2-y z_1)}
\end{equation*}
This gives us:
\begin{align*}
     \bm{R} = 
     \begin{pmatrix}
     1\\
     -\theta^i + \mathfrak{p}^2y +(\mathfrak{p} z_1 + \mathfrak{p}^2y_1 + \theta^i - \mathfrak{p}^2y) \tau\\
     \mathfrak{p} y + \mathfrak{p}^2z +(-\theta^j + \mathfrak{p}^2 z_1 - \mathfrak{p} y - \mathfrak{p}^2 z) \tau\\
     -\mathfrak{p} y + \mathfrak{p}^2y + (-\mathfrak{p} z_1 + \mathfrak{p}^2u_1 + \mathfrak{p} y - \mathfrak{p}^2 u) \tau
     \end{pmatrix}
\end{align*}
Let us set
\begin{align*}
    \tau' = -[\mathfrak{p}(\theta^{2j}y-\theta^{2i}z_1)+\mathfrak{p}^2(-y_1+z+y^2-y z_1)]=\\
    \mathfrak{p}(\theta^{2i}z_1 - \theta^{2j}y) + \mathfrak{p}^2 (y_1-z-y^2+y z_1)
\end{align*}
Then 
\begin{align}
\label{casea}
     \bm{R} = 
     \begin{pmatrix}
     \tau'\\
     (-\theta^i + \mathfrak{p}^2y) \tau' + \mathfrak{p} z_1 + \mathfrak{p}^2y_1 + \theta^i - \mathfrak{p}^2y \\
     (\mathfrak{p} y + \mathfrak{p}^2z) \tau' -\theta^j + \mathfrak{p}^2 z_1 - \mathfrak{p} y - \mathfrak{p}^2 z\\
     (-\mathfrak{p} y + \mathfrak{p}^2y) \tau' -\mathfrak{p} z_1 + \mathfrak{p}^2u_1 + \mathfrak{p} y - \mathfrak{p}^2 u 
     \end{pmatrix}
\end{align}
This finishes the proof of our statement.
\subsection{\texorpdfstring{Case 3: $\bm{P_1} = \bm{P_2} \mod \mathfrak{p}$}%
{Case 2: P_1 \neq P_2 \mod \mathfrak{p}}}
\begin{proposition}
\label{propap2}
Let ${\bm P_1} , {\bm P_2}, {\bm P_3}  \in V(K)$  be collinear and ${\bm P_1} \equiv_i {\bm P_2}$ where $1 \leq i \leq 3$ , then ${\bm P_1} \equiv_i {\bm P_3} $
\end{proposition}
\textbf{Proof}\\
If $\bm P_1 = \bm P_2$ then Proposition \ref{propap2} follows from Proposition \ref{propap}. Therefore we assume that
$\bm P_1 \neq \bm P_2$. 
We can follow  the proof in Proposition \ref{propap} until the equation (\ref{pointap1}).
Then we define:
\begin{equation*}
\label{pointap3}
\mathcal{T}{\bm P_2} = \left (\tilde{y} + \mathfrak{p}^j\tilde{\alpha},  \tilde{z} + \mathfrak{p}^j\tilde{\beta}, \tilde{u} + \mathfrak{p}^j\tilde{\gamma}  \right)
\end{equation*}
\begin{equation}
\label{pointap4}
\mathcal{T}{\bm P_3} = \left (\tilde{y} + \mathfrak{p}^j\tilde{\alpha}\tau,  \tilde{z} + \mathfrak{p}^j\tilde{\beta}\tau, \tilde{u} + \mathfrak{p}^j\tilde{\gamma} \tau \right)
\end{equation}
for some $( \tilde{\alpha},  \tilde{\beta},  \tilde{\gamma}) \in K^3$, $\tau \in k$, $j = i-1$ and
\begin{equation}
\label{ordprod}
\nu(\tilde \alpha)\nu(\tilde \beta)\nu(\tilde \gamma) = 0
\end{equation}
 since ${\bm P_1} \equiv_i {\bm P_2}$ and ${\bm P_1}, {\bm P_2}, {\bm P_3}$
are collinear.
Substituting (\ref{pointap4}) into (\ref{eq}) we get the following equalities:
\begin{equation*}
\label{eqap1}
-\theta^2(\tilde{y} + \mathfrak{p}^j \tilde{\alpha}t) + \theta^2 \mathfrak{p} (\tilde{y} + \mathfrak{p}^j \tilde{\alpha}t)^2 +  (\tilde{y} + \mathfrak{p}^j \tilde{\alpha}t)^3 +  (\tilde{z} + \mathfrak{p}^j \tilde{\beta}t)^3 + \theta  (\tilde{u} + \mathfrak{p}^j \tilde{\gamma} t)^3 = 0
\end{equation*}
\begin{align*}
-\theta^2\mathfrak{p}^j \tilde{\alpha} +\theta^2  \mathfrak{p}^{j+1}2\tilde{y}\tilde{\alpha} + 3\mathfrak{p}^{j}\tilde{y}^2\tilde{\alpha} +  3\mathfrak{p}^{j}\tilde{z}^2\tilde{\beta} + \theta 3\mathfrak{p}^{j}\tilde{u}^2\tilde{\gamma} + \\
(\theta^2  \mathfrak{p}^{2j+1}\tilde{\alpha}^2 + 3\mathfrak{p}^{2j}\tilde{y}\tilde{\alpha}^2 +  \mathfrak{p}^{2j}\tilde{z}\tilde{\beta}^2 + \theta 3\mathfrak{p}^{2j}\tilde{u}\tilde{\gamma}^2)t + \\
\mathfrak{p}^{3j}( \tilde{\alpha}^3 +  \tilde{\beta}^3 + \theta  \tilde{\gamma} ^3)t^2 = 0
\end{align*}
This can be represented as:
\begin{equation*}
l_0 + l_1\tau +l_2\tau^2=0
\end{equation*}
where
\begin{align*}
l_0 = -\theta^2 \tilde{\alpha} + 2\theta^2 \mathfrak{p}\tilde{y}\tilde{\alpha} -\theta^2 \mathfrak{p}^2 \tilde{y}^2\tilde{\alpha} -\theta^2 \mathfrak{p}^2 \tilde{z}^2\tilde{\beta} - \mathfrak{p}^2\tilde{u}^2\tilde{\gamma}\\
l_1 = \theta^2  \mathfrak{p}^{j+1}\tilde{\alpha}^2 -\theta^2 \mathfrak{p}^{j+2}\tilde{y}\tilde{\alpha}^2 -\theta^2 \mathfrak{p}^{j+2}\tilde{z}\tilde{\beta}^2 -\mathfrak{p}^{j+2}\tilde{u}\tilde{\gamma}^2\\
l_2 = \mathfrak{p}^{2j}( \tilde{\alpha}^3 +  \tilde{\beta}^3 + \theta  \tilde{\gamma} ^3)
\end{align*}
Since
\begin{equation}
\label{m}
l_0 + l_1 +l_2=0
\end{equation}
we get
\begin{equation*}
\tau = \frac{l_0}{l_2} = -\frac{l_1+l_2}{l_2} = -1 - \frac{l_1}{l_2}
\end{equation*}
\begin{equation}
\frac{l_1}{l_2}  = \frac{ \theta^2  \tilde{\alpha}^2 -\mathfrak{p} \theta^2  \tilde{y}\tilde{\alpha}^2 -\mathfrak{p} \theta^2  \tilde{z}\bar{\beta}^2 - \mathfrak{p} \tilde{u}\bar{\gamma}^2}{\mathfrak{p}^{j-1}( \tilde{\alpha}^3 +  \bar{\beta}^3 + \theta  \bar{\gamma} ^3)}
\end{equation}
We will be using (\ref{z3u3}) and it consequence  (\ref{sumzu}).\\
We have $\nu(\tilde \beta^3 + \theta \gamma^3) \leq 1$ and therefore
we consider the following  sub-cases.\\
\textbf{Sub-case i}
\begin{equation}
\label{case1}
\nu(\tilde \beta^3 + \theta \gamma^3) = 0
\end{equation}
If $j=0$ then we proceed as in the following. First, we have $-\theta^2 \tilde{\alpha}  +  \tilde{\alpha}^3 = 0 \mod \mathfrak{p}$. This contradicts to (\ref{case1}) because of (\ref{m})).
If $j > 0$ then
$\nu(\tilde \alpha) > 0$ and
\begin{equation}
\label{flow}
\nu \left (\frac{l_1}{l_2} \right ) \geq - (j - 2) \geq 0
\end{equation}
i.e.  $\nu(\tau) \geq 0$
and therefore
\begin{align}
\label{yzu}
\tilde{y} \equiv_j \tilde{y} + \mathfrak{p}^j\tilde{\alpha}\tau,\\
\label{yzu1}
\tilde{z} \equiv_j  \tilde{z} + \mathfrak{p}^j\tilde{\beta}\tau, \\
\label{yzu2}
\tilde{u} \equiv_j \tilde{u} + \mathfrak{p}^j\tilde{\gamma} \tau
\end{align}
Therefore
$\mathcal{T}{\bm P_1} \equiv_j \mathcal{T}{\bm P_3}$
and finally $\bm P_1 \equiv_i \bm P_3$.\\
\textbf{Sub-case ii}
\begin{equation}
\label{case2}
\nu(\tilde \beta^3 + \theta \gamma^3) = 1
\end{equation}
If $\nu(\tilde \alpha) = 0$ then $j=0$ and $\nu\left (\frac{l_1}{l_2} \right) = 1$ or $\nu(\tau) = 0$ proving the Sub-case ii.
Otherwise, if $\nu(\tilde \alpha) > 0$ then we can set $\tilde \alpha = \mathfrak{p} \hat \alpha, \; \hat \alpha \in K$
and we can represent
\begin{equation*}
\frac{l_1}{l_2}  = \frac{ \mathfrak{p} \theta^2  \hat{\alpha}^2 -\mathfrak{p}^2 \theta^2  \tilde{y}\hat{\alpha}^2 - \theta^2  \tilde{z}\bar{\beta}^2 -  \tilde{u}\bar{\gamma}^2}{\mathfrak{p}^{j-2}( \mathfrak{p}^3 \hat{\alpha}^3 +  \bar{\beta}^3 + \theta  \bar{\gamma} ^3)}
\end{equation*}
The equality  (\ref{case2}) and assumptions $\nu(\tilde \alpha) > 0$, (\ref{ordprod}) imply
\begin{equation*}
\tilde \beta = -\tilde \gamma \mod \mathfrak{p}
\end{equation*}
and
\begin{equation*}
\nu(\tilde \beta) = \nu(\tilde \gamma)  = 0
\end{equation*}
This implies:
\begin{equation*}
\frac{l_1}{l_2}  = \frac{\mathfrak{p} \theta^2  \hat{\alpha}^2 -\mathfrak{p}^2 \theta^2  \tilde{y}\hat{\alpha}^2 -(\theta^2  \tilde{z}+\tilde{u})\bar{\gamma}^2 + \mathfrak{p} h}{\mathfrak{p}^{j-2}( p^3 \hat{\alpha}^3 +  \bar{\beta}^3 + \theta  \bar{\gamma} ^3)}
\end{equation*}
for some $h \in K$.
Next, $\theta^2  \tilde{z}+\tilde{u} = \tilde{z}+\tilde{u} + \mathfrak{p} \theta^2 \tilde{z}$.
(\ref{sumzu}) and (\ref{case2}) implies that
\begin{equation*}
\frac{l_1}{l_2}  = \frac{\theta^2  \hat{\alpha}^2 -\mathfrak{p} \theta^2  \tilde{y}\hat{\alpha}^2 -h_1\bar{\gamma}^2 +  h}{\mathfrak{p}^{j-2}( \mathfrak{p}^2 \hat{\alpha}^3 +  h_2)}
\end{equation*}
where $h_1 = \frac{\theta^2  \tilde{z}+\tilde{u}}{\mathfrak{p}} \in K$ and $h_2 = \frac{\bar{\beta}^3 + \theta  \bar{\gamma} ^3}{\mathfrak{p}} \in K$, $\nu(h_2) = 0$
and then again we can have a flow like after (\ref{flow}).\\
Equations (\ref{yzu},\ref{yzu1},\ref{yzu2}) also show that if there is another collinear triple of points $\bm{P_1'}, \bm{P_2'}, \bm{P_3'} \in V(k)$, $\bm{P_1'} \neq \bm{P_2'}$ such that $\bm{P_1}\equiv_3 \bm{P_1'}$, $\bm{P_2}\equiv_3 \bm{P_2'}$ then $\bm{P_3}\equiv_3 \bm{P_3'}$.
If $\bm{P_1'}\equiv_3 \bm{P_2'}$ then $\bm{P_1}\equiv_3 \bm{P_2} \equiv_3 \bm{P_3}$ and by Proposition \ref{propap} $\bm{P_1'}\equiv_3 \bm{P_3'}$ completing the proof of admissibility of modulo $\mathfrak{p}^3$ on $V(k)$ in  the Case 3. This completes the proof of Theorem \ref{maintheorem}.
\section{Conjecture on cubic surfaces over number fields}
\label{conjecture}
Let $V \subset \mathbb{P}^3_L$ be a smooth cubic surface defined over a number field $L$ and such that $V(L)$ is non-empty.
Let $\nu_1$ and $\nu_2$ be two different places of $L$ such that the following holds.
\begin{enumerate}
    \item There exists a non-associative Moufang loop of point classes on $V(L_{\nu_1})$ (This can happen only if the residue field of $L_{\nu_1}$ has characteristic 3 \cite{dimitrikanevsky3})
    \item Let $\mathcal{O}$ be the ring of integers of $L_{\nu_2}$, which has a maximal ideal generated by uniformizing element $\pi$. Let $V_{\nu_2}$ be the Zariski closure of $V$ in $\mathbb{P}^3_{L_{\nu_2}}$. Let $\tilde V_{\nu_2}$ be the special fibre of $V_{\nu_2}$
    (this is a reduction of $V$ modulo $\pi$). Assume that $\tilde V_{\nu_2}$ is a cone over smooth cubic curve (over $\mathcal{O} \mod \pi$)
\end{enumerate}
Under these conditions we have
\begin{conj}
\label{conj}
There exists an admissible equivalence $\mathcal{A}$ on $V(L)$ such 
that $CML \; V(L)/\mathcal{A}$ is not associative.
\end{conj}
This conjecture is motivated by the fact that if a cubic surface over $\mathfrak{p}$-adic field $k$ has a ``bad'' reduction to a cone over a smooth cubic curve then generically the local evaluation map  $V(k) \rightarrow \mathbb{Q}/\mathbb{Z}$ takes as many values as possible (\cite{bright}). This can be used to show that for every point $\bm{P_1} \in V(L_{\nu_1})$ we can find a point $\bm{P_2} \in V(L_{\nu_2})$ such that the sum of local invariants over a certain adele (that includes $\bm{P_1}$ and $\bm{P_2}$) vanishes.
Using this observation one can expect (for example, in the spirit of Collioth-Th\'el\'ene's and Sansuc's conjecture H3 in the section 3.8 (\cite{clt})) that $V(L)$ is dense in $V(L_{\nu_1})$ and therefore a non-associative Moufang loop of point classes on $V(L_{\nu_1})$ induces a
non-associative Moufang loop of point classes over $V(L)$.

\end{document}